\documentclass{article}
\usepackage[utf8]{inputenc}
\usepackage{amsmath}
\usepackage{amsthm}
\usepackage{amssymb}
\usepackage{cite}
\usepackage{tikz}
\usepackage{xcolor}
\usepackage{ dsfont }
\usepackage{mathtools}

\newtheorem{example}[equation]{Example}
\newtheorem{theorem}[equation]{Theorem}
\newtheorem{definition}[equation]{Definition}
\newtheorem{lemma}[equation]{Lemma}
\newtheorem{corollary}[equation]{Corollary}

\newtheorem{observation}[equation]{Observation}
\newtheorem{remark}[equation]{Remark}
\newcommand\E{\mathbb{E}}

\DeclareMathOperator{\tr}{tr}
\DeclareMathOperator{\height}{ht}

\title{Asymptotics of local face distributions and the face distribution of the complete graph}
\author{Jesse Campion Loth\\[1mm]
Department of Mathematics\\
Simon Fraser University\\
Burnaby, BC, V5A 1S6, Canada}

\begin{document}
\maketitle

\begin{abstract}
    We are interested in the distribution of the number of faces across all the $2-$cell embeddings of a graph, which is equivalent to the distribution of genus by Euler's formula.  In order to study this distribution, we consider the local distribution of faces at a single vertex.  We show an asymptotic uniformity on this local face distribution which holds for any graph with large vertex degrees.  
    
    We use this to study the usual face distribution of the complete graph.  We show that in this case, the local face distribution determines the face distribution for almost all of the whole graph.  We use this result to show that a portion of the complete graph of size $(1-o(1))|K_n|$ has the same face distribution as the set of all permutations, up to parity.  Along the way, we prove new character bounds and an asymptotic uniformity on conjugacy class products.
\\

    
\end{abstract}

\section{Introduction}

\subsection{Background}

We are interested in $2-$cell embeddings of graphs on orientable surfaces, up to orientation preserving homeomorphism.  These embeddings are well-known to be in bijection with local rotations, which are cyclic orderings of the half-edges at each vertex in a graph.  The genus of an embedding is defined as the genus of the surface it is embedded on.  A classical problem is to determine the minimum and maximum genus over all embeddings of a graph.  For general graphs, this problem is very difficult: Thomassen\cite{thomassen1989graph} showed that determining the minimum genus of a graph is NP-complete.  There are some fixed classes of graphs for which the minimum and maximum genus is known, but there will be many embeddings with a genus between these two values.  Therefore, given a graph $G$, and a positive integer $g$, one could ask for the number of embeddings of $G$ of genus $g$.  This question has only been answered for families of graphs with especially nice structure, see \cite{gross2018calculating} and the references therein.  We will work with the number of faces instead of genus, as it better illustrates our results.  These two parameters determine each other through Euler's formula: A graph with $n$ vertices, $m$ edges, and $f$ faces, that is embedded on a surface of genus $g$, satisfies $n-m+f = 2 -2g$.

One focus of this paper will be the embeddings of the complete graph $K_n$.  It was shown by Nordhaus, Stewart and White \cite{nordhaus1971maximum} that there are always embeddings of $K_n$ with $1$ or $2$ faces. The minimum genus of $K_n$ is intimately connected with the Ringel-Youngs conjecture on map colouring.  This conjecture was proven by a string of several papers across the 1950s and 1960s, culminating in the final proof by Ringel and Youngs \cite{ringel1968solution}.  In our context, their proof gives that there are embeddings of $K_n$ with $\lfloor \tfrac{n(n-1)}{3} \rfloor$ faces for all $n \geq 3$.  However, the distribution of the number of faces across all embeddings of $K_n$ is still unknown, and there is no known fast method to compute this distribution in general.  A straightforward computation is blocked by the large number of embeddings: Even for small values of $n$ it proves to be very difficult, as $K_n$ has $(n-2)!^n$ embeddings.  The distribution of the number of faces across all embeddings of $K_7$ is calculated in \cite{beyer2016practical} and explicitly given in \cite{gross2018calculating}.  Due to the complexity of the problem, the distribution is not known for $n \geq 8$, and it is perhaps unrealistic to expect a `nice' formula for the number of embeddings of $K_n$ with $f$ faces.

Therefore, an approach previously used for this problem is a probabilistic one, which was termed \emph{random topological graph theory} by White \cite{white1994introduction}.  Stahl \cite{stahl1995average} gave an upper bound for the average number of faces across all the embeddings of $K_n$, and this was recently improved in \cite{loth2022randomcomplete}.  Another approach is to look at the local distribution of faces around a single vertex in the graph, across all of its embeddings.  This was introduced independently as the local genus polynomial/distribution by Gross et al. \cite{gross2016combinatorial} and by Chen and Reidys \cite{chen2016local}.  It was noted by F{\'e}ray \cite{feray2015combinatorial} that this local genus polynomial can be studied by using conjugacy class products.  This fact is used in \cite{loth2021random} to estimate the average number of faces of various graphs.

It is the purpose of this paper to study the asymptotics of local and non-local face distributions of graphs.  We will show that under certain conditions, the local face distribution satisifes an asymptotic uniformity.  We use this to show that almost all of the complete graph has a non-local face distribution satisfying an asymptotic uniformity.  We start by giving some basic results and definitions on random embeddings of graphs, before stating our main results.


\subsection{Random embeddings of graphs}
For a finite set $X$, we write $S_X$ for the symmetric group on $X$, and write $U_X$ for the uniform distribution on $X$.  A combinatorial map is a triple $m = (D,R,E)$, where $D$ is a set of even cardinality, and $R,E \in S_D$ such that $E$ is a fixed point free involution.  We say that a map is connected if $\langle R,E \rangle$ acts transitively on $[n]$; all the maps we consider in this paper will be connected.  We will be studying $2$-cell embeddings of graphs on oriented surfaces, up to orientation-preserving homeomorphism.  It is well known that any connected map gives an embedding of a graph on a surface.  In this context, the darts in $D$ are the half-edges of the graph and each cycle of $R$ is a cyclic ordering of the half-edges at a vertex.  The permutation $E$ then describes how the half-edges are joined together to make the edges of the graph.  We call $R \cdot E$ the \emph{face permutation} of $m$, since the cycles of $R \cdot E$ give the \emph{faces} of the embedding.  We start by defining how we can describe all the embeddings of a graph as maps.

Let $G$ be a simple graph with vertices $v_1, \dots, v_n$ and edges $E(G)$.  Let $D(G) := D_1(G) \cup D_2(G) \cup \dots \cup D_n(G)$ be a set of darts, with $| D_i(G) | = \deg(v_i)$ for each $i$.  The darts in $D_i(G)$ will correspond to the half-edges at vertex $v_i$ in $G$.  We define two sets of permutations $\mathcal{R}(G), \mathcal{E}(G)$ as follows:
\begin{itemize}
    \item $\mathcal{R}(G)$ is the set of all \emph{local rotations}.  This is the set of all permutations in $S_D$ with $n$ cycles, labelled $\pi_1, \dots, \pi_n$, with each cycle $\pi_i$ on the symbols $D_i(G)$.
    \item $\mathcal{E}(G)$ is the set of all \emph{edge schemes}.  This is the set of all fixed point free involutions in $S_D$, where for each $v_iv_j \in E(G)$ we have exactly one $2-$cycle containing a dart from $D_i(G)$ and a dart from $D_j(G)$.
\end{itemize}

We suppress notation and omit the $G$ when it is clear from the context.  It is well known that if we fix any edge scheme $E \in \mathcal{E}$, then the set $\{ (D,R,E) : R \in \mathcal{R} \}$ is in bijection with the set of all embeddings of $G$ given by local rotations, see \cite[$\mathsection$3.2.4]{mohar2001graphs} for details.  Therefore, each possible embedding of $G$ appears $| \mathcal{E} |$ times in $\{ (D,R,E) : R \in \mathcal{R}, E \in \mathcal{E} \}$.  We choose to study this set in which we vary the edge scheme instead of fixing it, as it will prove convenient in our analysis.  This explains the following observation:

\begin{observation}
    Let $\mathcal{M}(G) := \{ (D,R,E) : R \in \mathcal{R}(G), E \in \mathcal{E}(G) \}$.  
    Choosing an $m \in \mathcal{M}(G)$ uniformly at random gives a randomly chosen 2-cell embedding of $G$.
\end{observation}

We therefore define a \emph{random embedding} of $G$ as $m \sim U_{\mathcal{M}(G)}$, where $U_{\mathcal{M}(G)}$ is the uniform distribution on $\mathcal{M}(G)$.  Let $F_G = F_G(m)$ be the random variable on $m \sim U_{\mathcal{M}(G)}$ that is given by the number of faces in $m$.

\begin{figure}
    \centering
    \includegraphics[scale=0.7]{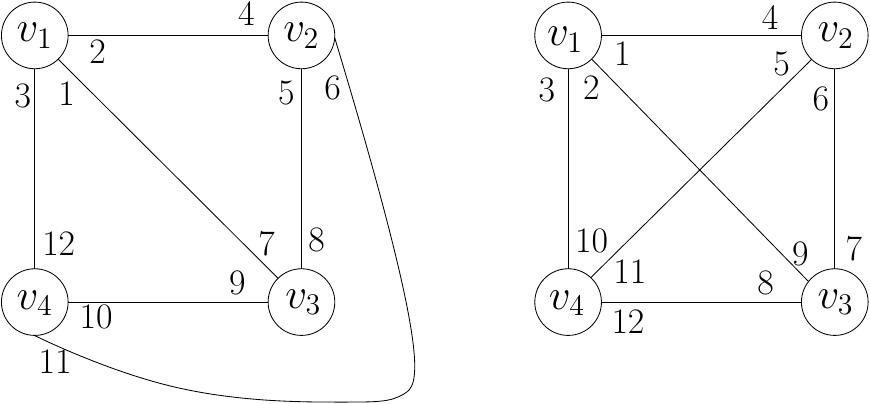}
    \caption{Two embeddings of $K_4$, drawn with clockwise local rotations.}
    \label{fig:firstexample}
\end{figure}

\begin{example}
    We give two examples of embeddings of $K_4$, the complete graph on $4$ vertices, in Figure \ref{fig:firstexample}.  The local rotations are drawn in clockwise order around each vertex, and $R \cdot E$ traces around the faces counter-clockwise.  The example on the left is a planar embedding.  The example on the right is an embedding on the torus, but we represent it with a drawing on the plane with edges crossing.  Both embeddings are given by maps where the set of darts $D$ is:
    $$
        D = D_1 \cup D_2 \cup D_3 \cup D_4 = \{1,2,3\} \cup \{4,5,6\} \cup \{7,8,9\} \cup \{10,11,12\}.
    $$
    The planar embedding has
    \begin{align*}
        R &= (2\,1\,3)(4\,6\,5)(7\,8\,9)(10\,11\,12), \\
        E &= (1\,7)(2\,4)(3\,12)(5\,8)(6\,11)(9\,10), \\
        R \cdot E &= (1\,12\,9)(2\,7\,5)(3\,4\,11)(6\,8\,10).
    \end{align*}
    The embedding on the torus has
     \begin{align*}
        R &= (1\,2\,3)(4\,6\,5)(7\,8\,9)(10\,11\,12), \\
        E &= (1\,4)(2\,9)(3\,10)(5\,11)(6\,7)(8\,12), \\
        R \cdot E &= (1\,9\,6\,11\,8\,2\,10\,5)(3\,4\,7\,12).
    \end{align*}
\end{example}

Fix any set $X$ and $Y \subseteq X$ and write some permutation $\alpha \in S_X$ in disjoint cycle notation.  Define the \emph{induced permutation} of $\alpha$ on $Y$ by deleting all symbols not contained in $Y$ from $\alpha$ in its presentation as a product of cycles, then removing all empty cycles.  We give a small example.  Let $Y = \{1,3,4,6\} \subset [7]$, and consider the permutation $\alpha := (1 \, 3 \, 2)(4 \, 6)(5 \, 7) \in S_7$ written in cycle notation.  Then the induced permutation of $\alpha$ on $Y$ is obtained by deleting the symbols $\{2,5,7\}$ from $\alpha$.  Since the cycle $(5 \, 7)$ only contains symbols in $\{2,5,7\}$, we remove this cycle altogether.  This gives the induced permutation $(1 \, 3)(4 \, 6) \in S_Y$.

Now fix some vertex $v_i$ in $G$ and write $v=v_i$ and $D_v = D_i$.

\begin{definition}
    Fix a simple graph $G$ and a vertex $v$, and let $m \in \mathcal{M}(G)$.  Let $\omega_v(m) \in S_{D_v}$ be the induced permutation of $R \cdot E$ on $D_v$.   The \emph{local face distribution} is the probability distribution on $S_{D_v}$ given by $\omega_v(m)$ for $m \sim U_{\mathcal{M}(G)}$.  We denote this probability distribution by $L_v(G)$.
\end{definition}  

When a map $m \in \mathcal{M}$ is fixed, we write $\omega_v = \omega_v(m)$. 

\begin{example}
    We continue our previous example on the two embeddings of $K_4$ given in Figure \ref{fig:firstexample}.  Let $v=v_1$.  The planar embedding has face permutation $(1\,12\,9)(2\,7\,5)(3\,4\,11)(6\,8\,10)$ and $D_v = \{1,2,3\}$.  It therefore has $\omega_v = (1)(2)(3)$.  Intuitively, this means the three faces incident with $v$ are distinct.  The toroidal embedding has face permutation $(1\,9\,6\,11\,8\,2\,10\,5)(3\,4\,7\,12)$ and $D_v = \{1,2,3\}$, so it has $\omega_v = (1 \, 2)(3)$.  This means that there are two faces incident with $v$, and that one of these faces visits $v$ twice.
\end{example}

We note that our definition of the local face distribution is very similar to the local genus polynomial introduced by Gross et al. \cite{gross2016combinatorial}, and the local genus distribution introduced by Chen and Reidys \cite{chen2016local}.  The difference is that their versions fix the value of $R - \pi_i$, whereas we let it run over all possibilities.

\subsection{Main Results}


We will start with a theorem on products of conjugacy classes in the symmetric group, which is of independent interest.  Let $A_X$ be the set of all even permutations in $S_X$, which is the alternating group, and let $A_X^c$ be the set of all odd permutations in $S_X$.  When $X=\{1,2,\dots,n\}$, write $S_n = S_X$ and $A_n = A_X, A_n^c = A_X^c$.  Let $C_\lambda$ be the set of all permutations of cycle type $\lambda = (\lambda_1, \dots, \lambda_m)$.  It is well known that the conjugacy classes of the symmetric group are given by $C_\lambda$ for all partitions $\lambda$ of $n$.  We say that a partition $\lambda$ is even (odd) if every permutation in $C_\lambda$ is even (odd). 

Let $\pi \sim U_{C_{\alpha}}, \sigma \sim U_{C_{\beta}}$, and $\gamma = \pi \sigma$.  Then let $P_{\alpha,\beta}$ be the probability distribution on $S_n$ given by $\gamma$.  The \emph{total variation distance} of two probability distributions $P,Q$ on a set $X$ is defined as $||P-Q|| = \tfrac{1}{2} \sum_{x \in X} |P(x) - Q(x)|$.

\begin{theorem} \label{thm:asymptoticunif}
    Let $\lambda$ be a partition of $n$ with $\ell$ parts of size $1$ such that $(n)$ and $\lambda$ have the same parity.  Then we have
    $$
        || P_{(n), \lambda} - U_{A_n} || \leq \frac{\max\{2, \ell + 1\}}{2\sqrt{n-1}}.
    $$
    The same statement holds with $U_{A_n^c}$ when $(n)$ and $\lambda$ have opposite parity.
\end{theorem}

We use Theorem \ref{thm:asymptoticunif} to study the local face distribution of general graphs.  We show that for any graph with a small average number of faces across all of its embeddings, the local face distribution at a vertex with a large degree is close to uniform.  For a set $X$, define the distribution $U_{X,p}$ to be the probability distribution on $S_X$ which takes value $\tfrac{2p}{|X|!}$ at odd permutations and $\tfrac{2(1-p)}{|X|!}$ at even permutations.  In other words, $U_{X,p}$ is uniform on $A_X$ and uniform on $A_X^c$, but is not necessarily uniform on all of $S_X$.  Let $G-\{v\}$ be the graph obtained from $G$ by removing vertex $v$ and all of the edges incident with $v$.

\begin{theorem} \label{thm:localunif}
    Let $v$ be a vertex of degree $d \geq 5$ in a graph $G$.  Set $p = \mathbb{P}[\omega_v(m) \text{ is odd}]$, where $m \sim U_{\mathcal{M}(G)}$.  Then we have
    $$
        || L_v - U_{D_v,p} ||^2 \leq \frac{2 \E[F_{G-\{v\}}]}{\sqrt{d-1}}.
    $$
\end{theorem}

In the case of the complete graph $K_n$, we show that for any vertex $v$ the portion of $K_n$ made up of the faces which contain at least one half-edge at $v$ has size $(1-o(1))|K_n|$.  These two results combine to give Theorem \ref{thm:knuniform}.  This says that a part of the complete graph of size $(1-o(1))|K_n|$ has a face distribution close to uniform.

Recall that the Stirling number of the first kind $c(n,k)$ is the number of permutations on $n$ symbols with $k$ cycles.  Define the probability distribution $P_{n,p}$ on $\{1,2,\dots ,n\}$ as
\begin{align*} 
    P_{n,p}(k) = \begin{cases}
        p \frac{2c(n,k)}{n!} : n + k \text{ odd}, \\
        (1-p) \frac{2c(n,k)}{n!} : n + k \text{ even}.
    \end{cases}
\end{align*}
In other words, $P_{n,p}$ is equal to the cycle distribution on $U_{[n],p}$.

\begin{theorem} \label{thm:knuniform}
    Set $m =(D,R,E) \sim U_{\mathcal{M}(K_n)}$, and let $D'(m) \subseteq D$ be the set of darts contained in faces incident with vertex $v$ in $m$.  Let $\alpha(m)$ be the induced permutation of $R \cdot E$ on $D'$, then we have:
    \begin{itemize}
        \item $\E[|D'|] = (1-o_n(1))|D|$.
        \item Let $Q_{n-1}$ be the distribution on $\{1,2,\dots,n-1\}$ for the number of cycles in $\alpha(m)$ and let $p = \mathbb{P}[\alpha(m) \text{ is odd}]$, then $|| Q_{n-1} - P_{n-1,p} || \rightarrow 0$.
    \end{itemize}
\end{theorem}

Our methods of proof are a mixture of probabilistic and algebraic ones.  We make use of recent bounds \cite{loth2022randomcomplete} on the average number of faces of a randomly chosen embedding of $K_n$, and techniques for calculating this average in general \cite{loth2021random}.  We also make use of some algebraic machinery: representations of the symmetric group.  We prove a new, simple bound for hook-shaped characters, and combine it with a general technique for calculating total variation distances between probability distributions \cite{chmutov2016surface, gamburd2006poisson}.  The representation theory is self-contained in Section \ref{sec:chartheory}.  In Sections \ref{sec:localclt} and \ref{sec:completegraphasymptotics} we use these representation theoretic results to study local and non-local face distributions.

\section{Character theory and class products} \label{sec:chartheory}

Recall that $C_\alpha$ is the conjugacy class containing all permutations of cycle type $\alpha$.  We abuse notation and write $C_X \subseteq S_X$ for the set of all full cycles on a set $X$.  If $\pi \sim U_{C_{\alpha}}, \sigma \sim U_{C_{\beta}}$, then recall that $P_{\alpha,\beta}$ is the probability distribution on $S_n$ given by $\gamma = \pi \sigma$.  The rest of the section is dedicated to the proof of Theorem \ref{thm:asymptoticunif}, and is the only part of the paper which uses character theory.

We follow a technique used by Gamburd \cite{gamburd2006poisson} and by Chmutov and Pittel \cite{chmutov2016surface}.  In their papers, they use a bound on the total variation distance due to Diaconis and Shahshahani \cite{diaconis1981generating}.  They combine it with a powerful, but asymptotic, bound on character values of the symmetric group due to Larsen and Shalev \cite{larsen2008characters}.  However, the character values in our application are simpler than those in the aforementioned papers, so we do not need to apply such general machinery.  It will also be convenient to have a more concrete character bound without asymptotic notation.  We will therefore start by estimating the character values we need with a simple bound.

Let $\chi^{\lambda}(\alpha)$ denote the character of the irreducible representation of $S_n$ indexed by $\lambda$ evaluated at a permutation of cycle type $\alpha$.  Write $f^\lambda$ for the dimension of the irreducible representation indexed by $\lambda$.

 \begin{lemma} \label{lem:hookshapecharbound}
    Let $\lambda$ be a partition of $n \geq 3$ with $\ell$ parts of size $1$, and suppose $1 \leq k \leq n-2$.  Then
    \begin{align*}
    \frac{|\chi^{(n-k,1^k)}(\lambda)|}{f^{(n-k,1^k)}} \leq \frac{\max\{2, \ell + 1\}}{n-1}.
    \end{align*}
\end{lemma}

\begin{proof}
We define a \emph{Murnaghan--Nakayama placement} of $\lambda$ into $(n-k,1^k)$ as a placement of $\lambda_i$ copies of $i$ into the the hook shape tableaux $(n-k,1^k)$, so that each row and column is weakly increasing, and for each $i$ all the boxes containing $i$ are connected.  The height $\height(T)$ of a placement $T$ is defined as the sum over all $i$, of one less than the number of columns of the connected section containing $i$.
The Murnaghan--Nakayama rule \cite[$\mathsection$7.17.3]{StanEC2} states that
$$
    \chi^{(n-k,1^k)}(\lambda) = \sum_T (-1)^{\height(T)},
$$
where the sum is over all Murnaghan--Nakayama placements $T$ of $\lambda$ into $(n-k,1^k)$.  Using this, the dimension of the hook shape $f^{(n-k,1^k)} = \chi^{(n-k,1^k)}(1^n)$ is easy to calculate: it is the number of ways of placing $\{1,2,\dots,n\}$ into the hook shape $(n-k,1^k)$ such that the first row and column is increasing.  The symbol $1$ must be placed in the top left box, then once we pick a subset of $k$ symbols to go into the first and only column of the shape, the rest of the placements are determined.  Since each of these placements has height $0$, we have $f^{(n-k,1^k)} =  \binom{n-1}{k}$.

Now write $\lambda = (\lambda_1, \lambda_2, \dots, \lambda_{t})$, where $\lambda_1 \geq \lambda_2 \geq \dots \geq \lambda_t$.  Since $\lambda$ has $\ell$ parts of size one, $\lambda_i \geq 2$ for $1 \leq i \leq t-\ell$.  A composition is defined as an ordered sequence of positive integers summing to $n$.  For any composition $\lambda'$ obtained by reordering the parts of $\lambda$, the Murnaghan--Nakayama rule tells us that $\chi^{(n-k,1^k)}(\lambda) = \chi^{(n-k,1^k)}(\lambda')$.  In order to make our calculations more convenient, we define a composition $\lambda' = (\lambda_t, \lambda_1, \dots, \lambda_{t-1})$.  Let $g^{(n-k,1^k)}(\lambda)$ be the total number of Murnaghan--Nakayama placements of the composition $\lambda'$ into $(n-k,1^k)$.  Then we have
$$|\chi^{(n-k,1^k)}(\lambda)| = |\chi^{(n-k,1^k)}(\lambda')| \leq g^{(n-k,1^k)}(\lambda).$$
We will spend the rest of this proof estimating
$$\frac{g^{(n-k,1^k)}(\lambda)}{\binom{n-1}{k}}.$$ 

We now perform a case analysis, based on the value of $k$ and the number of fixed points in $\lambda$.

Case 1: Suppose $\lambda$ has at most one part of size one.  Suppose $k \leq n/2$, the case for $k > n/2$ follows by symmetry since $g^{(n-k,1^k)}(\lambda) = g^{(k+1,1^{n-k-1})}(\lambda)$.  

There must be a copy of $1$ in the top left corner of the tableaux, then $r$ copies of $1$ directly below this, where $0 \leq r \leq \lambda_t-1$.  Every way of placing the other numbers into the tableaux is uniquely given by a subset of $\{\lambda_1, \dots, \lambda_{t-1} \}$ that sums to $k-r$, because these parts correspond to the symbols that are placed in the first and only column of the tableaux in weakly increasing order, then the remaining symbols are placed into the first and only row in weakly increasing order.  Note that the cardinality of such a subset is at most $(k-r)/2$, since $\lambda_i \geq 2$ for $1 \leq i \leq t-1$.  Suppose we have such a subset with
$$
    \lambda_{i_1} + \lambda_{i_2} + \dots + \lambda_{i_j} = k-r,
$$
where we could have $j=0$ and $k=r$.
There is at most one such subset for $k=1$, which is the empty subset, so $g^{(n-1,1)}(\lambda) \leq 1$ and the bound holds in this case.  Therefore, assume $k>1$.  Then we have that
\begin{align*}
    n &= \lambda_1 + \lambda_2 + \lambda_3 + \dots + \lambda_t \geq \lambda_t + (k-r) + 2(t - 1 - j) \\
    &\geq \lambda_t + (k-\lambda_1 + 1) + 2(t - 1 - j) \geq k - 1 + 2(t-j)
\end{align*}
Rearranging gives $j \geq k/2 - n/2 + t -1$, so this is a lower bound on the cardinality of a subset of the parts of $\lambda$ that sums to $k-r$ for some $0 \leq r \leq \lambda_t - 1$.  We can therefore bound $g^{(n-k,1^k)}(\lambda)$ by all subsets of $(\lambda_1, \dots, \lambda_{t-1})$ of cardinality between $k/2 - (n/2 - t + 1)$ and $k/2$.  We note that we could have that $k/2 - n/2 + t - 1 < 0$, but the following bounds still hold in this case.  Since $j$, the cardinality of the subset, must be an integer, we have
$$
    j \geq \lceil k/2 - (n/2 - t + 1) \rceil \geq \lfloor k/2 \rfloor - (\lfloor n/2 \rfloor - t + 1).
$$
Letting $s = \lfloor n/2 \rfloor - t + 1 \geq 0$, we obtain 
\begin{align*}
    g^{(n-k,1^k)}(\lambda)  \leq \sum_{j=\lfloor k/2 \rfloor-s}^{\lfloor k/2\rfloor} \binom{t-1}{j} =  \sum_{j=0}^s \binom{\lfloor n/2 \rfloor - s}{\lfloor k/2 \rfloor -j}\leq \sum_{j=0}^s \binom{\lfloor n/2 \rfloor - s}{\lfloor k/2 \rfloor-j} \binom{s}{j}.
\end{align*}

If $\lfloor k/2 \rfloor < j < s$, then $\binom{\lfloor n/2 \rfloor - s}{\lfloor k/2 \rfloor -j} = 0$.  Therefore, we may upper bound the sum by changing the upper limit of summation to $\lfloor k/2 \rfloor$.  Using Vandermonde's identity, we have
\begin{align*}
    g^{(n-k,1^k)}(\lambda)  &\leq \sum_{j=0}^{\lfloor k/2 \rfloor} \binom{\lfloor n/2 \rfloor - s}{\lfloor k/2 \rfloor-j} \binom{s}{j} = \binom{\lfloor n/2 \rfloor}{\lfloor k/2 \rfloor}.
\end{align*}

\begin{figure}
    \centering
    \includegraphics[scale=0.6]{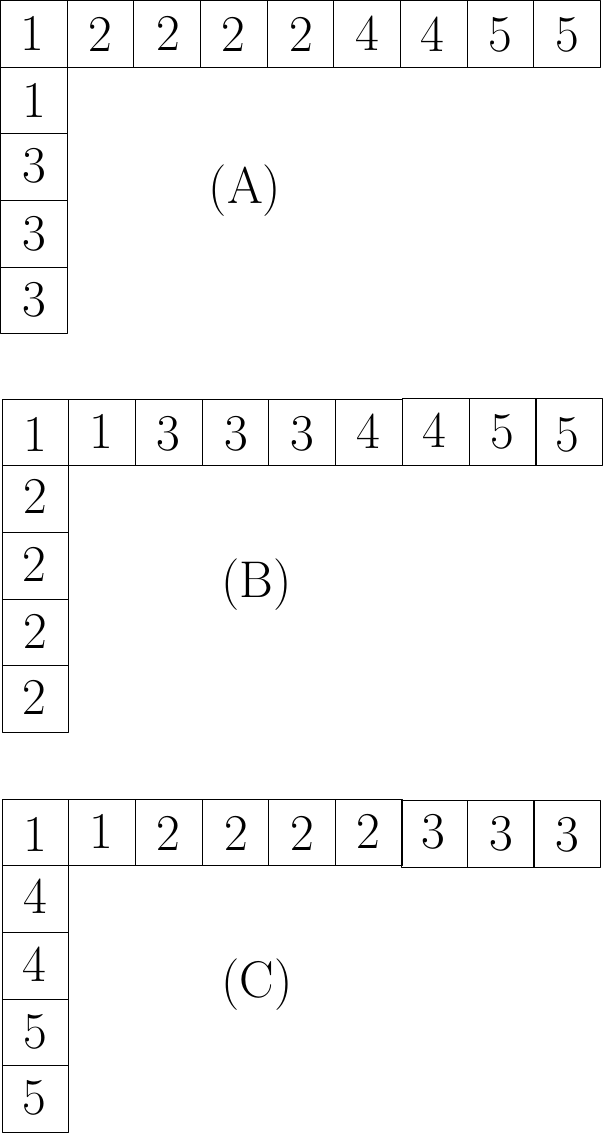}
    \caption{Let $n=13$ and $\lambda = (4,3,2,2,2)$, then $\lambda' = (2,4,3,2,2)$ and $g^{(9,1^4)}(\lambda)$ is equal to the number of Murnaghan--Nakayama placements of $\{1,1,2,2,2,2,3,3,3,4,4,5,5\}$ into the tableaux of shape $(9,1^4)$.  We give an example of Case 1 of the proof of Lemma \ref{lem:hookshapecharbound}.  In this case, we bound the number of Murnaghan--Nakayama placements by the number of subsets of $(\lambda_1, \lambda_2, \lambda_3, \lambda_4) = (4,3,2,2)$ that sum to some value in $\{3,4\}$.  The valid subsets in this case are $\{\lambda_2\}, \{\lambda_1\}, \{\lambda_3, \lambda_4\}$, and they refer to placements $(A),(B),(C)$ respectively.}
    \label{fig:tableauxplacement}
\end{figure}

Therefore, for $k$ even we have 
$$\frac{g^{(n-k,1^k)}(\lambda)}{f^{(n-k,1^k)}} \leq \frac{\binom{\lfloor n/2 \rfloor}{k/2}}{\binom{n-1}{k}}.$$  
For $k$ odd, since $k<n/2$, we have
$$
    \frac{g^{(n-k,1^k)}(\lambda)}{f^{(n-k,1^k)}} \leq \frac{\binom{\lfloor n/2 \rfloor}{(k-1)/2}}{\binom{n-1}{k}} \leq \frac{\binom{\lfloor n/2 \rfloor}{(k-1)/2}}{\binom{n-1}{k-1}}
$$

For $k$ even we obtain a bound of the following form for $i=k$, and $i \geq 2$.  For $k$ odd we have $k \geq 3$, and we obtain a bound of the following form for $i = k-1 \geq 2$.
\begin{align*}
    \frac{g^{(n-k,1^k)}(\lambda)}{f^{(n-k,1^k)}} \leq \frac{\binom{\lfloor n/2 \rfloor}{i/2}}{\binom{n-1}{i}}
    &\leq \frac{n(n-2)(n-4)\dots(n-i+2)}{i(i-2)(i-4)\dots2} \frac{i!}{(n-1)(n-2)\dots(n-i)}\\
    &= \frac{n}{n-i} \frac{i!}{i(i-2)(i-4)\dots 2} \frac{(n-2)(n-4)\dots(n-i+2)}{(n-1)(n-2)\dots(n-i+1)} \\
    &= \frac{n}{n-i} \frac{(i-1)(i-3)(i-5)\dots1}{(n-1)(n-3)(n-5)\dots(n-i+1)} \\
    &\leq \frac{n}{(n-i)(n-1)} \left( \frac{i-1}{n-3} \frac{i-3}{n-5} \dots \frac{3}{n-i+1} \right) \\
    &\leq \frac{2}{n-1} \leq \frac{\max\{2, \ell + 1\}}{n-1},
\end{align*}
where in the last line we've used the fact that $\tfrac{n}{n-i} \leq 2$.

Case 2: The partition $\lambda$ has $\ell \geq 2$ fixed points, and $k=1$ or $n-2$.
First suppose $k=1$, as the case $k=n-2$ gives the same bound by symmetry.  The single $1$ must be placed in the top left of the tableaux.  The square just below this $1$ in the tableaux must be filled with an $i$ with $\lambda'_i = 1$.  There are $\ell-1$ choices of such a part, then the rest of tableaux placements are determined.  This gives that $g^{(n-1,1)}(\lambda) = \ell-1$.  In fact since the height of each of these tableaux is $0$, we have that $\chi^{(n-1,1)}(\lambda) = g^{(n-1,1)}(\lambda) = \ell-1$.  Therefore, for this case we have
$$
    \frac{g^{(n-1,1)}(\lambda)}{f^{(n-1,1)}} = \frac{\ell-1}{n-1} < \frac{\ell + 1}{n-1}.
$$

Case 3: The partition $\lambda$ has $\ell \geq 2$ fixed points and $2\leq k \leq n-3$.  
This covers the remaining possible cases for $\lambda$.  We proceed by induction on $n$ and $\ell$, using the previous two cases.  Suppose that for all $\mu \vdash n-1$ with $\ell -1$ fixed points and all $1 \leq k \leq n-2$ we have
$$
    \frac{g^{(n-k,1^k)}(\mu)}{f^{(n-k,1^k)}} \leq \frac{\ell}{n-2}.
$$

Firstly note that if $\ell=n$ then $g^{(n-k,1^k)}(\lambda) = g^{(n-k,1^k)}(1^n) = f^{(n-k,1^k)}$ and the bound is trivial, so we assume $\ell \leq n-2$. 
 Note that $\lambda_{t-1} = 1$, we write $\lambda - \lambda_{t-1}$ for the partition of $n-1$ obtained by removing $\lambda_{t-1}$ from $\lambda$.  The single copy of $t$ corresponding to $\lambda_{t-1}$ must be placed at the end of the first row or at the bottom of the first column in the tableaux.  If it is placed in the far right, then there are $g^{(n-k-1, 1^k)}(\lambda - \lambda_{t-1})$ ways to place the remaining symbols into the tableaux.  If this $m$ is placed in the bottom cell then there are $g^{(n-k, 1^{k-1})}(\lambda - \lambda_{t-1})$ ways to place the remaining symbols.  Note that since $f^{(n-k,1^k)} = \binom{n-1}{k}$, we have $f^{(n-k,1^k)} = f^{(n-k-1,1^k)} + f^{(n-k,1^{k-1})}$.  Therefore,
\begin{align*}
    \frac{g^{(n-k,1^k)}(\lambda)}{f^{(n-k,1^k)}} \leq \frac{g^{(n-k-1, 1^k)}(\lambda - \lambda_{t-1})+g^{(n-k, 1^{k-1})}(\lambda - \lambda_{t-1})}{f^{(n-k-1,1^k)} + f^{(n-k,1^{k-1})}}.
\end{align*} 
By our inductive assumption, and the fact that $\ell \geq 2$, we have
\begin{align*}
    \frac{g^{(n-k-1, 1^k)}(\lambda - \lambda_{t-1})}{f^{(n-k-1,1^k)}}, \frac{g^{(n-k, 1^{k-1})}(\lambda - \lambda_{t-1})}{f^{(n-k,1^{k-1})}} \leq \frac{\ell}{n-2}.
\end{align*}
Therefore, since $\ell \leq n-2$ we have
\[
    \frac{g^{(n-k,1^k)}(\lambda)}{f^{(n-k,1^k)}} \leq \frac{\ell}{n-2} \leq \frac{\ell+1}{n-1}. \qedhere
\]
\end{proof}

We are now ready to prove Theorem \ref{thm:asymptoticunif}.

\begin{proof}[Proof of Theorem \ref{thm:asymptoticunif}]

We prove the formula in the case that $(n), \lambda$ are partitions of the same parity.  The same arguments hold for partitions of opposite parity, with $U_{A_n}$ replaced by $U_{A_n^c}$.

Chmutov and Pittel \cite{chmutov2016surface} give the proof for $\alpha = 2^{n/2}$, but their proof generalises.  We present this proof, closely following their structure and adding a few extra details.  Suppose $\alpha, \beta \vdash n$ have the same parity, and write $P = P_{\alpha,\beta}$.  Note that $P$ is a class function that takes the value zero at any odd permutation. 

    The starting point is the Cauchy--Schwarz inequality:
    \begin{align*}
        \left( \sum_{x \in X} \vert Q(x)  T(x) \vert \right)^2 \leq \left( \sum_{x \in X} |Q(x)|^2 \right) \left( \sum_{x \in X} |T(x)|^2  \right).
    \end{align*}
    Letting $X = S_n$, $Q = P - U_{A_n}$ and $T = 1$ in the above formula gives
    \begin{align*}
        || P - U_{A_n} ||^2 &= \left( \frac{1}{2} \sum_{\omega \in S_n} \vert P(\omega) - U_{A_n}(\omega) \vert \right)^2 \\
        &\leq \frac{|S_n|}{4}  \sum_{\omega \in S_n} |P(\omega) - U_{A_n}(\omega)|^2.
    \end{align*}

    We can evaluate this sum using the Plancherel Theorem \cite[Theorem $\mathsection$15.2(2)]{terras1999fourier}.  Let $\rho^{\lambda}$ be the irreducible representation of $S_n$ indexed by $\lambda$, and recall that $f^\lambda$ is the dimension of this representation.  For a function $P:S_n \rightarrow \mathbb{C}$ the Fourier transform of $P$ at a representation $\rho$ is defined as $\hat{P}(\rho) = \sum_{\omega \in S_n} P(\omega) \rho(\omega)$.  Then we have
    \begin{align*}
        |S_n| \sum_{\omega \in S_n} |P(\omega) - U_{A_n}(\omega)|^2 = \sum_{\lambda \vdash n} f^\lambda tr[(\hat{P}(\rho^\lambda) - \hat{U}_{A_n}(\rho^\lambda))(\hat{P}(\rho^\lambda) - \hat{U}_{A_n}(\rho^\lambda))^*].
    \end{align*}
     We continue by evaluating $\hat{P}(\rho^\lambda) - \hat{U}_{A_n}(\rho^\lambda)$ for each $\lambda \vdash n$.
    
    The case $\lambda = (n)$ refers to the trivial representation, and $\lambda = (1^n)$ refers to the sign representation.  Therefore, $\rho^{(n)}(\omega) = 1$ and $\rho^{(1^n)}(\omega) = (-1)^{\text{sgn}(\omega)}$ for all $\omega \in S_n$.  By the definition of the Fourier transform, for $\lambda = (n)$ we have
    $$
        \hat{P}(\rho^{(n)}) - \hat{U}_{A_n}(\rho^{(n)}) = \sum_{\omega \in S_n} P(\omega) - \sum_{\omega \in A_n} \frac{1}{|A_n|} = 1-1 = 0.
    $$

    By similar reasoning, and using the fact that $(n)$ and $\lambda$ have the same parity, $\hat{P}(\rho^{1^n}) - \hat{U}_{A_n}(\rho^{1^n}) = 0$.  Now let $\lambda \neq (n), (1^n)$.  Since distinct irreducible characters are orthogonal, we have
    $$
         \sum_{\omega \in A_n} \chi^{\lambda}(\omega) + \sum_{\omega \in A_n^c} \chi^{\lambda}(\omega) = \sum_{\omega \in S_n} \chi^{\lambda}(\omega) = \sum_{\omega \in S_n} \chi^{\lambda}(\omega) \chi^{(n)}(\omega) = 0.
    $$

    Also, using the sign representation,
    $$
       \sum_{\omega \in A_n} \chi^{\lambda}(\omega) - \sum_{\omega \in A_n^c} \chi^{\lambda}(\omega) = \sum_{\omega \in S_n} \chi^{\lambda}(\omega) \chi^{(1^n)}(\omega) = 0.
    $$

    Since for any two complex numbers $a,b$, we have $a+b=a-b=0$ implies that $a=b=0$, we have for any $\lambda \neq (n), (1^n)$,
    $$
    tr\left( \hat{U}_{A_n}(\rho^\lambda) \right) =  \sum_{\omega \in A_n} \chi^{\lambda}(\omega) = 0.
    $$
   By an application of Schur's Lemma \cite[Lemma $\mathsection$9.1]{james2001representations} we have that $\sum_{\omega \in A_n} \rho^{\lambda}(\omega)$ is some constant times the identity matrix.  However, since the trace of this matrix is zero, this implies that $\hat{U}_{A_n}(\rho^\lambda) = 0$.
    
      Therefore, we arrive at the following bound for the total variation distance, often called the Diaconis--Shahshahani upper bound Lemma \cite{diaconis1981generating}.
    \begin{align} \label{eqn:totalvardistbound}
        || P - U_{A_n} ||^2 \leq \frac{1}{4} \sum_{\lambda \neq (n), (1^n)} f^\lambda \tr(\hat{P}(\rho^\lambda) \hat{P}(\rho^\lambda)^*).
    \end{align}
    
    We therefore continue by evaluating $\hat{P}(\rho^\lambda)$.  Firstly note that since the Fourier transform is multiplicative over convolutions of probability distributions, we have that $\hat{P} = \hat{U}_{C_\alpha} \cdot \hat{U}_{C_\beta}$.  It therefore suffices to calculate $\hat{U}_{C_\alpha}$ for $\alpha \vdash n$.  By definition we have
    \begin{align*}
        \hat{U}_{C_\alpha}(\rho^\lambda) &= \frac{1}{|C_\alpha|} \sum_{\omega \in C_\alpha} \rho^\lambda(\omega).
    \end{align*}
    By Schur's lemma the matrix on the right hand side is some constant times the identity matrix, and $\tr(\rho^\lambda(\omega)) = \chi^\lambda(\omega)$.  This gives
    $$
        \hat{U}_{C_\alpha}(\rho^\lambda) = \frac{1}{f^\lambda} \chi^\lambda(\alpha) \mathbb{I}_{f^\lambda},
    $$
    where $\mathbb{I}_{f^\lambda}$ is the identity matrix of dimension $f^\lambda$.  Therefore, we have
    $$
        \hat{P}(\rho^\lambda) = \frac{\chi^\lambda(\alpha)  \chi^\lambda(\beta) }{(f^\lambda)^2} \mathbb{I}_{f^\lambda}.
    $$
    Finally, since the irreducible characters of the symmetric group are real valued, we have
    \[
        \hat{P}(\rho^\lambda) \hat{P}(\rho^\lambda)^* = \left( \frac{\chi^\lambda(\alpha)  \chi^\lambda(\beta) }{(f^\lambda)^2} \right)^2 \mathbb{I}_{f^\lambda}.
    \]
    Taking the trace in the previous formula gives
    $$
        \tr\left( \hat{P}(\rho^\lambda) \hat{P}(\rho^\lambda)^* \right) = \left( \frac{\chi^\lambda(\alpha)  \chi^\lambda(\beta) }{(f^\lambda)^2} \right)^2 f^\lambda.
    $$
    Substituting this into \eqref{eqn:totalvardistbound} gives
$$
    || P_{\alpha, \beta} - U_{A_n} ||^2 \leq \frac{1}{4} \sum_{\lambda \neq (n), (1^n)} \left( \frac{\chi^{\lambda}(\alpha) \chi^{\lambda}(\beta)}{f^\lambda} \right)^2.
$$
We are interested in the case of $\alpha = (n)$.  It is only possible to place $n$ copies of $1$ into a tableaux of shape $\lambda$ so that they form a hook if $\lambda$ is itself a hook shape.  If $\lambda = (n-k, 1^k)$, there is just a single such placement, which is of height $k$.  Therefore, by the Murnaghan--Nakayama rule we have
$$
\chi^\lambda((n)) =  \begin{cases} (-1)^k \text{ if } \lambda = (n-k,1^k),\\
0 \text{ else}.
    \end{cases}
$$
This simplifies the formula considerably:
$$|| P_{(n), \lambda} - U_{A_n} ||^2 \leq \frac{1}{4} \sum_{k=1}^{n-2} \left( \frac{(-1)^k \chi^{(n-k,1^k)}(\lambda)}{f^{(n-k,1^k)}} \right)^2. $$
Using the bound from Lemma \ref{lem:hookshapecharbound} we obtain:
\[
    || P_{(n), \lambda} - U_{A_n} ||^2 \leq \frac{1}{4} \sum_{k=1}^{n-2} \left( \frac{\max\{2, \ell + 1\}}{n-1} \right)^2 \leq \frac{(\max\{2, \ell + 1\})^2}{4(n-1)}. \qedhere
\]
\end{proof}

\section{Asymptotic uniformity of local face distributions} \label{sec:localclt}

We continue by applying this conjugacy class product result to local face distributions.  We start with a few definitions.  Recall the definition of a random embedding of $G$ as $m=(D,R,E) \sim U_{\mathcal{M}(G)}$.  Fix some vertex $v = v_i$, and write $D_v = D_i$, $\pi_v = \pi_i$ and $d_v = d_i$.

Let $R- \pi_v$ be the permutation in $S_D$ obtained from $R$ by replacing the cycle $\pi_v$ with $d_v$ fixed points.  Let $\sigma_v(m) \in S_{D_v}$ be the induced permutation of $(R-\pi_v) \cdot E$ on $D_v$.  Suppose we obtain a new map from $m$ by splitting the vertex $v$ into $d_v$ vertices with one half-edge incident with each, as shown in Figure \ref{fig:addingvertextofirstexample}.  Then $(R-\pi_v) \cdot E$ is the permutation in $S_D$ which is obtained by walking counter-clockwise around the faces in this new map.   Recall that $\omega_v(m) \in S_{D_v}$ is the induced permutation of $R \cdot E$ on $D_v$.

\begin{example}
Consider $v=v_1$ in the embedding of $K_4$ given in Figure \ref{fig:addingvertextofirstexample}.   Here we have:
\begin{align*}
    \omega_v &= (1\,2)(3), \\
    (R-\pi_v)\cdot E &= (1)(2)(3)(4\,6\,5)(7\,8\,9)(10\,11\,12) \cdot (1\,4)(2\,9)(3\,10)(5\,11)(6\,7)(8\,12)\\
    &= (1\,4\,7\,12\,3\,10\,5)(2\,9\,6\,11\,8), \\
    \sigma_v &= (1\,3)(2).
\end{align*}

\begin{figure}
    \centering
    \includegraphics[scale=0.7]{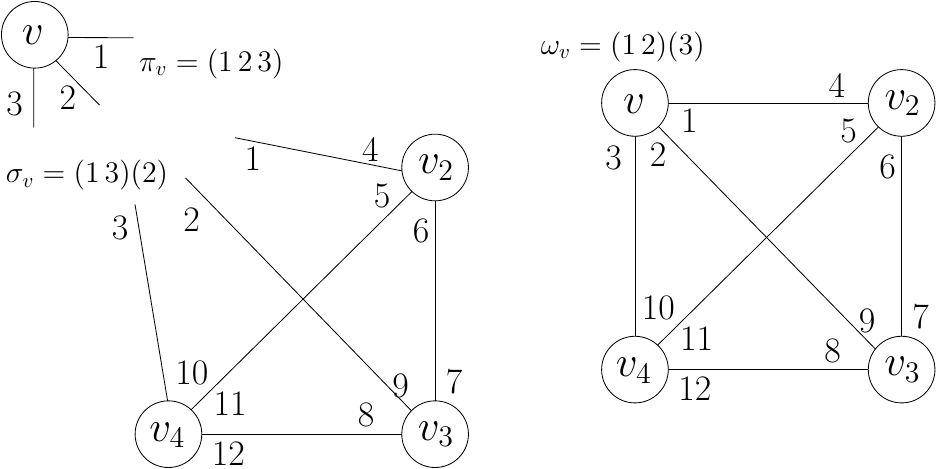}
    \caption{Continuing from Example 1, we give a pictorial explanation of why $\pi_v \cdot \sigma_v = \omega_v$.}
    \label{fig:addingvertextofirstexample}
\end{figure}

Notice that we have
$$
    \pi_v \cdot \sigma_v = (1\,2\,3) \cdot (1\,3)(2) = (1\,2)(3) = \omega_v.
$$
We will see that this correspondence between these three permutations is true in general.
\end{example}

We work in the symmetric group algebra $\mathbb{C} S_{D_v}$ on the set of darts in $D_v$.  Let $K_v(G) = \sum_{m \in \mathcal{M}(G)} \omega_v(m) \in \mathbb{C} S_{D_v}$, and let $K_v = K_v(G)$ when the graph is clear from the context.  For a partition $\lambda$ of $[d_v]$, let $C_\lambda(D_v)$ be the set of permutations in $S_{D_v}$ of cycle type $\lambda$.  Then we write $K_\lambda = \sum_{\tau \in C_\lambda(D_v)} \tau \in \mathbb{C} S_{D_v}$ for the \emph{class sum} corresponding to $\lambda \vdash [d_v]$.  We start by showing how the permutations $\omega_v, \sigma_v$ and $\pi_v$ are related.  Some of this correspondence is also outlined in \cite{loth2021random} using different notation.
    
\begin{lemma}  \label{lem:claim1}
    $$K_v(G) = \frac{1}{(d_v-1)!} K_{(d_v)} \sum_{m \in \mathcal{M}(G)} \sigma_v(m).$$
\end{lemma}

    \begin{proof}
      Take some map $m=(D,R,E) \in \mathcal{M}(G)$.  By definition, $\sigma_v(m)$ is only dependent on $R-\pi_v$ and $E$.  In particular, for every other possible full cycle $\pi_v'$ there is some map $m'$ with local rotation $(R-\pi_v) \cdot \pi_v'$ and $\sigma_v(m') = \sigma_v(m)$.
    
    It is therefore enough to show that for any fixed map $m \in \mathcal{M}(G)$, we have that $\omega_v(m) = \pi_v \sigma_v(m)$.  First we have
    $$
        R \cdot E = \pi_v \cdot (R-\pi_v) \cdot E.
    $$
    Therefore, the induced permutation of both sides on $D_v$ is the same. 
    
    The induced permutation of $R \cdot E$ on $D_v$ was defined as $\omega_v$.  The induced permutation of $(R-\pi_v) \cdot E$ is $\sigma_v$.  Since $\pi_v$, when viewed as a permutation in $S_D$, is just a single cycle only supported on $D_v$, the induced permutation of $\pi_v \cdot (R-\pi_v) \cdot E$ on $D_v$ is $\pi_v \sigma_v$.
    
    Therefore, $\pi_v \sigma_v = \omega_v$ as required. 
    \end{proof}

We continue with an observation on the set of all possible $\sigma_v(m)$ across all $m \in \mathcal{M}(G)$.  We write $Z(\mathbb{C} S_{D_v})$ for the centre of the group algebra $\mathbb{C} S_{D_v}$.  

\begin{lemma}  \label{lem:claim2}
    $$\sum_{m \in \mathcal{M}(G)} \sigma_v(m) \in Z(\mathbb{C} S_{D_v}).$$
\end{lemma} 
\begin{proof}
    It is enough to show $\sum_{m \in \mathcal{M}(G)} \sigma_v(m) = \tau \sum_{m \in \mathcal{M}(G)} \sigma_v(m) \tau^{-1}$ for any $\tau \in S_{D_v}$.  This is equivalent to showing that there is a bijection from $\mathcal{M}(G)$ to itself, taking $m$ to $m'$, so that the map $m'$ satisfies $\sigma_v(m') = \tau \sigma_v(m) \tau^{-1}$.  Fix some embedding $m=(D,R,E)\in \mathcal{M}(G)$.  Then since $R-\pi_v$ has all the symbols in $D_v$ as fixed points, we have
    $$
        (R-\pi_v) \cdot \tau E \tau^{-1}  = \tau (R-\pi_v) \cdot E \tau^{-1}.
    $$
    Now notice that conjugating a permutation by some $\tau \in S_n$ commutes with taking an induced permutation.  That is, suppose for some $\alpha \in S_n$ the induced permutation of $\alpha$ on some $Y \subseteq [n]$ is $\beta$.  Then the induced permutation of $\tau \alpha \tau^{-1}$ on $\tau (Y)$ is $\tau \beta \tau^{-1}$.
    
    Therefore, the embedding given by $m' = (D,R,\tau E \tau^{-1}) \in \mathcal{M}(G)$ has $\sigma_v(m') = \tau \sigma_v(m) \tau^{-1}$.  The map $m$ can be obtained by simply applying the same process with $\tau$ and $\tau^{-1}$ swapped, so this gives the required bijection. 
\end{proof}

Combining the previous two Lemmas gives us a key result.

\begin{lemma} \label{lem:sigmaomegarelation}
    Fix a graph $G$ and a vertex $v$.  Then there exist constants $a_{\lambda,v}$ for all $\lambda \vdash n$ such that
    \begin{align*}
        K_{v} = \sum_{\lambda \vdash [d_v]} a_{\lambda,v} K_{(d_v)} K_\lambda
    \end{align*}
\end{lemma}
\begin{proof}
 It is well known that $Z(\mathbb{C}S_{D_v})$ is spanned by the class sums $K_\lambda$ for $\lambda \vdash [d_v]$, see \cite[Proposition $\mathsection$12.22]{james2001representations}.  Therefore, the group algebra element $\sum_{m \in \mathcal{M}(G)} \sigma_v(m)$ can be expressed as
\begin{align*}
    \sum_{m \in \mathcal{M}(G)} \sigma_v(m) &= \sum_{\lambda \vdash d} b_{\lambda,v} K_\lambda.
\end{align*}
for some constants $b_{\lambda,v}$.  Using Lemmas \ref{lem:claim1} and \ref{lem:claim2} gives
\begin{align*}
    K_v(G) &= K_{(d_v)} \frac{1}{(d_v-1)!} \sum_{m \in \mathcal{M}(G)} \sigma_v(m) \\
    &=K_{(d_v)} \sum_{\lambda \vdash d} a_{\lambda,v} K_\lambda =  \sum_{\lambda \vdash d} a_{\lambda,v} K_{(d_v)} K_\lambda
\end{align*}
for some constants $a_{\lambda,v}$, proving the result.
\end{proof}

Therefore, in order to study the local face distribution using $K_v(G)$, we can use our results on conjugacy class products.  Recall that our upper bounds on $||P_{(d),\lambda} - U_{A_d}||$ and $||P_{(d),\lambda} - U_{A_d^c}||$ from Theorem \ref{thm:asymptoticunif} are only good when the number of parts of size $1$ in $\lambda$ is relatively small.  By the previous discussion, in our application each permutation in $C_\lambda$ in the product $C_{(d)} C_\lambda$ will be $\sigma_v(m)$ for some $m \in \mathcal{M}$.  We therefore need to show that the number of fixed points in $\sigma_v$ is, on average, relatively small.  We begin by bounding the number of fixed points in $\sigma_v$ by the number of faces in an embedding of a smaller graph.  For $m \sim U_{\mathcal{M}(G)}$, let $P_v$ be the random variable for the number of fixed points in $\sigma_v(m)$, and recall that $F_G(m)$ is the random variable for the number of faces in $m$.  Let $G-\{v\}$ be the graph obtained by removing vertex $v$ and all of the edges incident with $v$ from $G$.

\begin{lemma} \label{lem:fpbound}
     For any vertex $v$ we have that $\E[P_v] \leq \E[F_{G - \{v\}}]$.
\end{lemma}
\begin{proof}
    Let $m \sim U_{\mathcal{M}(G)}$.  There are $d_v$ 2-cycles in $E$ which contain a symbol in $D_v$.  Remove all the $2d_v$ symbols contained in these 2-cycles from $D$ to obtain $D'$.  Let $R'$ and $E'$ be the induced permutations of $R$ and $E$ respectively on the set $D'$.  We call the map obtained in this way $m-\{v\} = (D',R',E')$.  By definition this process gives $m-\{v\} \sim U_{\mathcal{M}(G - \{v\})}$.
    
    Recall that $\sigma_v$ is the induced permutation of $(R - \pi_v)\cdot E$ on $D_v$.  Therefore, the number of fixed points in $\sigma_v$ is certainly at most the number of cycles in $(R-\pi_v)\cdot E$.  By our construction of $m-\{v\} = (D',R',E')$, the face permutation $R' \cdot E'$ has the same number of cycles as $(R-\pi_v) \cdot E$.  Therefore, the permutation $(R-\pi_v)\cdot E$ has $F(m-\{v\})$ cycles, meaning that $P_v(m) \leq F(m-\{v\})$.  Taking expectation gives $\E[P_v] \leq \E[F_{G-\{v\}}]$, as required.
\end{proof}

We are ready to prove the first main theorem of the section.  Recall that $U_{D_v,p}$ is the probability distribution on $S_{D_v}$ which takes value $\tfrac{2p}{n!}$ at odd permutations in $S_{D_v}$ and $\tfrac{2(1-p)}{n!}$ at even permutations in $S_{D_v}$.  In other words, $U_{D_v,p}$ is uniform on $A_{D_v}$ and uniform on $A_{D_v}^c$, but is not necessarily uniform on all of $S_{D_v}$.  Recall that $L_v$ is the local face distribution of $G$ at vertex $v$, given by taking $\omega_v(m)$ for $m \sim U_{\mathcal{M}(G)}$.
 
\begin{proof}[Proof of Theorem \ref{thm:localunif}]
    By Markov's inequality \cite[Lemma $\mathsection$7.2(7)]{grimmett2020probability} we have $\mathbb{P}[P_v \geq x \E[P_v]] \leq \frac{1}{x}$ for any $x > 0$.  By Lemma \ref{lem:fpbound} we therefore have that $\mathbb{P}[P_v \geq x \E[F_{G-\{v\}}]] \leq \frac{1}{x}$.  Let $x^2 = 2\sqrt{d-1}/\E[F_{G-\{v\}]}]$, then since $d \geq 5$ we have that $x \E[F_{G-\{v\}}] \geq 2$.
    
    Therefore, by Theorem \ref{thm:asymptoticunif}, if a partition $\lambda$ with the same parity as $(d_v)$ has at most $x \E[F_{G-\{v\}}]$ parts of size $1$ then 
    $$||P_{(d),\lambda} - U_{A_d}|| \leq \frac{\max\{2, x\E[F_{G-\{v\}}]\} }{2\sqrt{d-1}} = \frac{x\E[F_{G-\{v\}}]}{2\sqrt{d-1}}.$$
    The same result holds for partitions $\lambda$ with opposite parity to $(d_v)$, with $U_{A_d}$ replaced by $U_{A_d^c}$.
    
    Then using Lemma \ref{lem:sigmaomegarelation} we have
    \begin{align*}
        || L_v - U_{d,p} || &\leq \sum_{\lambda \text{ even}} \mathbb{P}[\sigma_v(m) \in C_\lambda] ||P_{(d), \lambda} - U_{A_n}|| \\
        &+ \sum_{\lambda \text{ odd}} \mathbb{P}[\sigma_v(m) \in C_\lambda] ||P_{(d), \lambda} - U_{A_n^c}|| \\
        &\leq \mathbb{P}[P_v < x \E[F_{G-\{v\}}]] \frac{x\E[F_{G-\{v\}}]}{2\sqrt{d-1}} + \mathbb{P}[P_v \geq x \E[F_{G-\{v\}}]] \\
        &\leq \frac{x\E[F_{G-\{v\}}]}{2\sqrt{d-1}} + \frac{1}{x}.
    \end{align*}

    Substituting in $x^2 = 2\sqrt{d-1}/\E[F_{G-\{v\}]}]$ gives the result.
\end{proof}

In particular, Theorem \ref{thm:localunif} says that if $G$ has a large degree vertex $v$, and $\E[F_{G-\{v\}}]$ is small, then its local face distribution at $v$ will be close to a distribution derived from a mixture of the uniform distributions on even and odd permutations, with some unknown mixing parameter $p = \mathbb{P}[\omega_v(m) \text{ is odd}]$.  Note that although the face permutation $R \cdot E$ has the same parity for every $m \in \mathcal{M}$, the permutation $\sigma_v(m)$ varies in parity across different $m \in \mathcal{M}$.  This means the mixing parameter $p = \mathbb{P}[\omega_v(m) \text{ is odd}]$ in the statement of Theorem \ref{thm:localunif} is difficult to calculate, as it is hard to tell when $\sigma_v(m)$ will be odd or even.

\section{Asymptotics of complete graph embeddings} \label{sec:completegraphasymptotics}

We give a concrete example of Theorem \ref{thm:localunif}.  We use a recent result that bounds the average number of faces in the complete graph.

\begin{theorem} \cite[Theorem 1.4]{loth2022randomcomplete} \label{thm:knavgfaces}
    For sufficiently large $n$, we have that $\E[F(K_n)] \leq 4\log(n)$.
\end{theorem}

This allows us to give a special case of Theorem \ref{thm:localunif} for the complete graph.

\begin{corollary} \label{cor:knlocaluniform}
    Let $v$ be a vertex in the complete graph $K_n$.  For $m \sim U_{\mathcal{M}(G)}$, let $p = \mathbb{P}[\omega_v(m) \text{ is odd}]$.  Then for sufficiently large $n$ we have
    \begin{align*}
        || L_v - U_{D_v,p} ||^2 &\leq \frac{8\ln(n-1)+2}{\sqrt{n-2}} \\
        &\rightarrow 0.
    \end{align*}
\end{corollary}
\begin{proof}
    Combining Theorem \ref{thm:localunif} and Theorem \ref{thm:knavgfaces} gives the result.
\end{proof}

This means that the local face distribution of the complete graph at any vertex is asymptotically uniform on even permutations, and on odd permutations, but not necessarily on the whole symmetric group.  These results give us an understanding of the local behaviour of complete graph embeddings at a vertex.  We continue by analysing the behaviour of these embeddings across separate vertices.

For some dart $e \in D$, write $F_e = F_e(m)$ for the face containing $e$ in the embedding $m = (D,R,E)$.  This means that $F_e$ is a cycle of $R \cdot E$, containing some number of darts in $D$.  We say $v \in F_e$ if there is at least one dart $d \in D_v$ that is contained in the face $F_e$, and $v \notin F_e$ otherwise.  Write $|F_e(m)|$ for the number of darts contained in the face $F_e(m)$.

\begin{example}
    We give an example of these definitions on the left embedding of $K_4$ given in Figure \ref{fig:firstexample}.  This is a planar embedding with face permutation $R \cdot E = (1\,12\,9)(2\,7\,5)(3\,4\,11)(6\,8\,10)$.  We have $F_2 = (2 \, 7 \, 5)$, and $F_2 = F_5 = F_7$.  Since $2,7$ and $5$ are darts incident with vertices $v_1,v_3$ and $v_2$ respectively, we have $v_1 \in F_2$, $v_3 \in F_2$ and $v_2 \in F_2$.  Since these are the only darts in the face $F_2$, we have that $v_4 \notin D_2$.
\end{example}

We first give a way of expressing $\E[F_{K_{n}}]$ as a sum in terms of $|F_e(m)|$.

\begin{lemma}  \label{lem:averagefacesreformulation}
    Let $m \sim U_{\mathcal{M}(G)}$.  Then for any fixed dart $e \in D(K_n)$ we have
    \begin{align*}
    \E[F_{K_{n}}] = n(n-1) \sum_{k=1}^{n(n-1)} \frac{1}{k} \mathbb{P}[|F_e(m)| = k].
    \end{align*}
\end{lemma}
\begin{proof}
    Recall that $F(m)$ is the number of faces in $m \in \mathcal{M}(K_n)$.  For each dart $e \in D(K_n)$ and $m \in \mathcal{M}(K_n)$, let $X_e(m)$ be the Bernoulli random variable taking value $1$ with probability $1/|F_e(m)|$.  Then we have
    \begin{align*}
        \sum_{e \in D} \E[X_e(m)] = \sum_{e \in D} 1/|F_e(m)| = F(m). 
    \end{align*}
    
Therefore, summing over all $m$ and using the symmetry of $e \in D$ gives
\begin{align*}
   \E[F_{K_n}] &= \frac{1}{|\mathcal{M}(K_n)|} \sum_{m \in \mathcal{M}(K_n)} \sum_{e \in D} \E[X_e(m)] \\
   &=  \sum_{e \in D} \frac{1}{|\mathcal{M}(K_n)|} \sum_{m \in \mathcal{M}(K_n)} \E[X_e(m)] \\
   &= \sum_{e \in D} \E[1/|F_e|] = n(n-1) \E[1/|F_e|]. \qedhere
\end{align*}
\end{proof}

\begin{theorem} \label{thm:vertexinfaceprob}
  Let $v$ be a vertex in $K_n$, let $e$ be a dart not in $D_v$, and write $\mathrm{e}$ for Euler's constant.  Then for $m \sim U_{\mathcal{M}(G)}$ we have the following for sufficiently large $n$:
  \begin{align*}
    \mathbb{P}[v \notin F_e] &\leq \frac{4\ln(n-1)}{\mathrm{e}(n-1)} \\
    &\rightarrow 0.
  \end{align*}
\end{theorem}
\begin{proof} 
    Given $m' \sim U_{M(K_{n-1})}$, we give a process to add a vertex $v=v_n$ to extend it to $m \sim U_{M(K_n)}$.  
    \begin{itemize}
        \item Take $m' \sim U_{M(K_{n-1})}$.  Suppose $D'_i$ is the set of darts incident with vertex $v_i$ for $i=1,2,\dots,n-1$, and that $R' = \pi'_1 \pi'_2 \dots \pi'_{n-1}$ with $\pi'_i \in S_{D'_i}$ for each $i$.
        \item For each $i=1,2,\dots,n-1$ let $D_i = D'_i \cup \{d'_i\}$, then place $d'_i$ after a symbol chosen uniformly at random from $D'_i$ in the cycle $\pi_i'$ to obtain $\pi_i$.
        \item Let $D_n = \{e_1,e_2, \dots, e_{n-1}\}$ be a set of darts of size $n-1$, and let $\pi_n$ be a uniformly random full cycle on $D_n$.
        \item Let $R = \pi_1 \pi_2 \dots \pi_n$.
        \item Pick a uniformly random bijection $g: [n-1] \rightarrow [n-1]$.  Then let \newline $E = E' \cdot (d'_1 \, e_{g(1)}) (d'_2 \, e_{g(2)}) \dots (d'_{n-1} \, e_{g(n-1)})$.
    \end{itemize}
      It is immediate that this process outputs $m \sim U_{M(K_n)}$.

      Fix an embedding $m'=(D',R',E') \in \mathcal{M}(K_{n-1})$ and let $e$ be a dart in $D'$.  Suppose that there are $f_i$ darts in $F_e(m')$ contained in $D_i'$ for $i=1,\dots,n-1$.  Write $d_e(m') := (f_1,\dots,f_{n-1})$; see Figure \ref{fig:addingvertexbreakingfaces} for an example of this sequence.  At the second step of this process there are $n-2$ total choices of place to add the new dart $d_i'$ into the cycle $\pi_i'$.  If at any of the vertices one of the $f_i$ placements just after a dart in $F_e(m')$ is chosen, then in the newly obtained map $m \in \mathcal{M}(K_n)$ we have that $v \in F_e(m)$.  Therefore, we have that $n \notin F_e(m)$ only if one of the other $n-2-f_i$ places is chosen to add the new dart at each vertex $v_i$. 
    This probability is only dependent on the numbers $f_i$, so we have
    $$
        \mathbb{P}[v \notin F_e(m) \;\big|\; d_e(m') = (f_1, \dots, f_{n-1})] \leq \prod_{i=1}^{n-1} \left(1-\frac{f_i}{n-2}\right).
    $$

\begin{figure}
        \centering
        \includegraphics[scale=0.7]{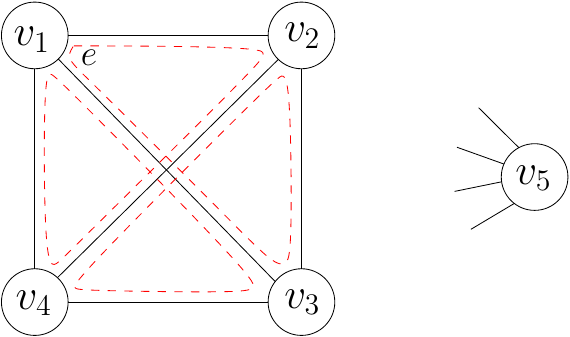}
        \caption[Adding a vertex to an embedding of $K_4$.]{This is an embedding of $K_4$ with two faces on the torus.  It has two faces, one of length $4$ and one of length $8$.  The face of length $8$ contains the dart $e$, and is outlined in red.  For the dart $e$ on the diagram of the map $m$, we have $d_e(m) = (2,2,2,2)$.  If we add $v_5$ to this embedding to give an embedding of $K_5$, there are $3$ places at each of $v_1,v_2,v_3,v_4$ to place the new half-edges.  Exactly one choice of placement at each vertex doesn't intersect $F_e$.  Therefore, in this case for $v=v_5$, $\mathbb{P}[v \notin F_e] = (1/3)^4$.}
        \label{fig:addingvertexbreakingfaces}
    \end{figure}
    
    Suppose $|F_e(m')| = \sum_{i=1}^{n-1} f_i = k$, then we have
    \begin{align*}
        \prod_{i=1}^{n-1} \left(1-\frac{f_i}{n-1}\right) &\leq \prod_{i=1}^{n-1} \mathrm{e}^{-f_i/(n-2)} \\
        &= \mathrm{e}^{-\sum_{i=1}^{n-1}f_i/(n-2)} =
        \mathrm{e}^{-k/(n-2)}.
    \end{align*}
    where we have used the bound $1-x \leq \mathrm{e}^{-x}$.
    
    Using total probability we obtain
    \begin{align}  \label{eq:halfwayprob}
        \mathbb{P}[v \notin F_e(m)] &= \sum_{k=1}^{(n-1)(n-2)} \mathbb{P}[|F_e(m')| = k] \mathbb{P}[v \notin F_e(m) \;\big|\; |F_e(m')| = k] \nonumber\\
        &\leq \sum_{k=1}^{(n-1)(n-2)} \mathbb{P}[|F_e(m')| = k] \mathrm{e}^{-k/(n-2)}.
    \end{align}   
    
    Let $m' \sim U_{\mathcal{M}(K_{n-1})}$.  
 Rearranging Lemma \ref{lem:averagefacesreformulation} and using Theorem \ref{thm:knavgfaces}, we obtain the following for sufficiently large $n$:
    \begin{align} \label{eq:log/n2thing}
        \sum_{k=1}^{(n-1)(n-2)} \frac{1}{k} \mathbb{P}[|F_e(m')| = k] \leq \frac{4\ln(n-1)}{(n-1)(n-2)} .
    \end{align}
    
    Then for $m \sim U_{\mathcal{M}(K_n)}$, using Equations \ref{eq:halfwayprob} and \ref{eq:log/n2thing} we have
    \begin{align*}
        \mathbb{P}[v \notin F_e(m)] &\leq \sum_{k=1}^{(n-1)(n-2)} \mathbb{P}[|F_e(m')| = k ] \mathrm{e}^{-k/(n-2)} \\
        &= \sum_{k=1}^{(n-1)(n-2)} \frac{1}{k} \mathbb{P}[|F_e(m')| = k] \left( k\mathrm{e}^{-k/(n-2)} \right) \\
        &\leq \sum_{k=1}^{(n-1)(n-2)} \frac{1}{k} \mathbb{P}[|F_e(m')| = k] \left( \frac{n-2}{\mathrm{e}} \right) \\
        &\leq \frac{4\ln(n-1)}{\mathrm{e}(n-1)} \rightarrow 0.
    \end{align*}
    The penultimate line uses the fact that $x\mathrm{e}^{-x/(n-2)}$ attains its maximum at $x=n-2$.
\end{proof}  

As a corollary to Theorem \ref{thm:vertexinfaceprob}, we can estimate the proportion of the embedding which is not changed when we alter the local rotation at a vertex $v$.

\begin{corollary} \label{cor:unchangedpartofgraph}
    Let $v$ be a vertex in $K_n$.  Let $D'(m) \subseteq D$ be the set of darts contained in faces incident with vertex $v$ in some $m \in \mathcal{M}(K_n)$.  Then $\E[|D'|] = (1-o_n(1))|D|$.
\end{corollary}
\begin{proof}
    We can decompose $|D'(m)|$ into indicator variables for each dart then use linearity of expectation to obtain
    $$
        \E[|D'(m)|] = \sum_{e \in D} \mathbb{P}[v \in F_e].
    $$
    The result then follows from Theorem \ref{thm:vertexinfaceprob}.
\end{proof}

We combine Theorem \ref{thm:localunif} with Corollary \ref{cor:unchangedpartofgraph} to obtain an asymptotic uniformity for embeddings of complete graphs.  Recall that the Stirling number of the first kind $c(n,k)$ is the number of permutations on $n$ symbols with $k$ cycles.  Recall that the probability distribution $P_{n,p}$ on $\{1,2,\dots ,n\}$ was defined as
\begin{align*} 
    P_{n,p}(k) = \begin{cases}
        p \frac{2c(n,k)}{n!} : n + k \text{ odd}, \\
        (1-p) \frac{2c(n,k)}{n!} : n + k \text{ even}.
    \end{cases}
\end{align*}
In other words, $P_{n,p}$ is equal to the cycle distribution on $U_{[n],p}$.

\begin{proof} [Proof of Theorem \ref{thm:knuniform}]
    Take $m \sim U_{M(K_n)}$ and fix a vertex $v$.  Recall that we defined $D'(m) \subseteq D$ as the set of all darts contained in faces that have at least one dart in $D_v$.  Let $\alpha(m)$ be the induced permutation of $R \cdot E$ on $D'$, then by construction we have that $c(\alpha(m))$ = $c(\omega_v(m))$ is just the number of faces incident with $v$.  Also, by Corollary \ref{cor:unchangedpartofgraph}, $\E[|D'|] = (1-o(1)) |D|$.   By Corollary~\ref{cor:knlocaluniform} the total variation distance between $L_{v}$ and $U_{D_v, p}$ tends to zero.  Using the triangle inequality, we have
    \begin{align*}
        ||L_v - U_{D_v, p}|| &= \frac{1}{2} \sum_{\omega \in S_{D_v}} |L_v(\omega) - U_{D_v, p}(\omega)| \\
        &= \frac{1}{2} \sum_{k=1}^n \sum_{\substack{\omega \in S_{D_v} \\ c(\omega) =k}} |L_v(\omega) - U_{D_v, p}(\omega)| \\
        &\geq \frac{1}{2} \sum_{k=1}^n \Big\vert \sum_{\substack{\omega \in S_{D_v} \\ c(\omega) =k}} \left(L_v(\omega) - U_{D_v, p}(\omega)\right)\Big\vert \\
        &= \frac{1}{2} \sum_{k=1}^n |Q_{n-1}(k) - P_{n-1,p}(k)|
    \end{align*}
    Therefore, we also have $||Q_{n-1} - P_{n-1,p}|| \rightarrow 0$, completing the proof.
\end{proof}

\begin{remark}
    Although $\E[|D'|] = (1-o(1))|D|$, we have that $\E[|D|-|D'|] \rightarrow \infty$.  There is no $D'$ so that Theorem \ref{thm:knuniform} holds with $\E[|D|-|D'|] \rightarrow 0$, since $K_n$ has embeddings with a quadratic number of faces.  Therefore, Theorem \ref{thm:knuniform} states that although there could be many faces, the distribution of longer faces satisfies an asymptotic uniformity.
\end{remark}

\section*{Acknowledgements}

The author would like to thank Kevin Halasz, Tomáš Masařík, Bojan Mohar and Robert Šámal for helpful discussions on the topic.  The author would also like to thank Amarpreet Rattan for carefully reading a draft of many of the proofs in the paper.

\bibliography{cites}
\bibliographystyle{plain}

\end{document}